\documentclass{svproc}
\usepackage[numbers]{natbib}

\usepackage{url}

\usepackage[utf8]{inputenc}
\usepackage{amsfonts, amsmath, amssymb, mathrsfs}
\usepackage{subfig}
\usepackage{mathtools}
\usepackage{geometry, graphicx}
\usepackage[hidelinks]{hyperref}
\usepackage{cleveref}
\usepackage{enumitem}
\hypersetup{colorlinks={true},linkcolor={blue}}
\usepackage{booktabs}
\usepackage{lipsum}
\usepackage{tikz-cd}
\usepackage{verbatim}
\usepackage{bm}

\usepackage{setspace}
\usepackage{indentfirst}

\usepackage{algorithm, algorithmic}

\crefformat{equation}{(#2#1#3)}

\graphicspath{ {./images/} }

\DeclareMathOperator*{\argmin}{arg\,min}
\DeclarePairedDelimiterX{\norm}[1]{\lVert}{\rVert}{#1}

\newcommand\set[1]{\left\{ #1 \right\}}

\begin{document}
\mainmatter              % start of a contribution
\title{Physics-informed neural networks via stochastic Hamiltonian dynamics learning}
\titlerunning{NeuralPMP}  % abbreviated title (for running head)
\author{Chandrajit Bajaj \inst{1}\inst{3}  \and Minh Nguyen\inst{2}\inst{3}}
\authorrunning{Minh, Bajaj}

\institute{Department of Computer Science \\
The University of Texas at Austin\\
\and
Department of Mathematics\\
The University of Texas at Austin\\
\and
Oden Institute for Computational Engineering and Sciences\\
The University of Texas at Austin
}

%\tocauthor{First Author, Second Author, Third Author and Fourth Author}
%
\maketitle              % typeset the title of the contribution

\begin{abstract}
In this paper, we propose novel learning frameworks to tackle optimal control problems by applying the Pontryagin maximum principle and then solving for a Hamiltonian dynamical system. Applying the Pontryagin maximum principle to the original optimal control problem shifts the learning focus to reduced Hamiltonian dynamics and corresponding adjoint variables. Then, the reduced Hamiltonian networks can be learned by going backwards in time and then minimizing loss function deduced from the Pontryagin maximum principle's conditions. The learning process is further improved by progressively learning a posterior distribution of the reduced Hamiltonians. This is achieved through utilizing a variational autoencoder which leads to more effective path exploration process. We apply our learning frameworks called NeuralPMP to various control tasks and obtain competitive results.
\end{abstract}
\section{Introduction}
Learning optimal control solution models for unknown dynamical systems, where an objective cost functional is optimized over space-time state-action sequences, is known to be computationally similar to convergently learning an optimal state-action policy in reinforcement learning. The actions on optimal value paths achieve maximum reward or minimum regret \cite{Todorov2006-jq}. The solutions to such reinforcement learning methods involve stable numerical solution techniques for high-dimensional and high variance constrained optimization. Several reinforcement learning algorithms have been devised attempting to compute the optimal control paths through guided search and policy learning \cite{abbeel2006using,deisenroth2011pilco,levine2014learning,heess2015learning,mnih2015human,dpo}. 

In this paper, we take a more comprehensive  approach. We rely mainly on two things: the physical nature or environment of the unknown dynamical system and the necessary requirements of the optimal control sequences or trajectories subjected to these physical dynamics \cite{boltyanskiy1962mathematical}. We learn the underlying reduced Hamiltonian dynamics and then utilize this data-dependent reduced Hamiltonians to optimally control the dynamical system for targeted objectives. We focus on the discrete-time with uneven time steps and continuous-time versions of the optimal control optimization problem. Through Pontryagin maximum principle, we reduce the optimal control problem to the data specific form of the Hamiltonian dynamics learning. The control is now encoded in the adjoint variable. Therefore, the new reinforcement learning task is to learn through this adjoint variable and the reduced Hamiltonian. This is fundamentally different from the optimal policy driven reinforcement learning, which focuses on deriving optimal action and value functions \cite{mnih2015human,gu2016continuous}. In particular, we first learn a correct reduced Hamiltonian of the controlled dynamics. This is obtained via the application of the Pontryagin maximum principle to the original optimal control problem. The reduced Hamiltonian deep network can be trained by going backwards in time and minimizing a Pontryagin maximum principle-based loss function. Moreover, the learning process is improved through variational inference and by learning a posterior distribution of the reduced Hamiltonians. This process leads to a filtration of the generalized coordinates during the neural network training.

\subsection{Related works}
Throughout this paper, we use PMP as the abbreviation for Pontryagin maximum principle, the fundamental principle for our main learning frameworks.\\

\textbf{Model-based reinforcement learning algorithms.}\quad Several model-based algorithms \cite{Wang2019-fk} have been developed for reinforcement learning. Unlike their model-free counterparts \cite{Mnih2016-cy, mnih2015human, TRPO, PPO}, model-based approaches model the outside environment and then use suitable mathematical formulations to facilitate the finding of agent's optimal actions. Two popular model-based approaches include the Dyna-style algorithms \cite{Sutton1990-qy} and algorithms focusing on policy search with backpropagation through time. While the first approach (Dyna-style algorithms) models the environment to generate imaginary data for actual policy training, the second approach, with notable methods such as PILCO \cite{deisenroth2011pilco}, computes the analytic gradient (closed-form formulas) of the reinforcement learning objective with respect to the policy via explicit modeling of the problem/environment.\\

\textbf{PMP-based framework on optimal control.}\quad Many works are developed based on the Pontryagin maximum principle (PMP) to tackle optimal control problems. For instance, Pontryagin Neural Networks (PoNNs) \cite{PoNNs} use PMP to transform the original problem into a Two-Point Boundary Value ordinary differential equation (ODE). PoNNs then use neural-network regression to learn such an ODE. On the other hand, Pontryagin Differentiable Programming \cite{PDP_end} and AI Pontryagin \cite{AI_pontryagin} use PMP to derive the exact analytic gradient of a certain loss function with respect to the parameters of the control to be learned. This method is similar to the second popular approach of the model-based reinforcement learning. We remark that all of these works leverage PMP only to derive closed-form formulas that can then be used to train the actual learning architecture. In these cases, PMP is not actively incorporated in the learning framework. \\

\textbf{PMP-based deep learning framework.}\quad There is also a line of research works \cite{Li2017-wy, zhang2019you} that use iterative PMP-based algorithms to train the Hamiltonian and then derive the optimal control variable. However, this line of works aims to improve the traditional supervised deep learning training rather than solving the optimal control tasks themselves. Their goals also focus on optimizing the parameters of deep learning architectures rather than the optimal paths of the control problem.\\

\textbf{Physics-informed dynamic learning.} \quad Several physics-informed learning architectures aim at learning hidden dynamics by incorporating physics biases. Deep Lagrangian Networks (DeLaN) \cite{lutter2019deep}, for example, injects the Lagrangian mechanics bias into the supervised learning process to improve the performance on difficult trajectories. Sympletic-ODE Net \cite{Zhong2019-df}, on the other hand, learns the unknown dynamics of a controlled system with Hamiltonian bias. Similar to our framework, Sympletic-ODE exploits the Hamiltonian dynamics for more accurate learning:
\begin{equation}
\begin{bmatrix}\dot{q}\\ \dot{p}\end{bmatrix} = f_{\theta}(q, p, u) = \begin{bmatrix}\frac{\partial H_{\theta_1, \theta_2}}{\partial p}\\ -\frac{\partial H_{\theta_1, \theta_2}}{\partial q}\end{bmatrix} + \begin{bmatrix} 0 \\ g_{\theta_3}(q) \end{bmatrix}u
\end{equation}
where 
\begin{equation}
H_{\theta_1, \theta_2}(q,p) = \frac{1}{2}p^TM^{-1}_{\theta_1}(q)p+V_{\theta_2}(q)
\end{equation}
However, these prior works focus more on the forward dynamics including the controlled dynamics rather than finding a control that optimizes a functional cost. Furthermore, in these works, the truth dynamics/trajectories are provided and the learnings are often conducted in a supervised manner.\\

\textbf{Neural ODE.} \quad Similar to our work, Sympletic-ODE Net builds on the Neural ODE architecture \cite{chen2018neuralode}, which proposes to model a neural network by an ordinary differential equation and learns the rate of change $\frac{dy}{dx} = f(x, y, \theta)$ rather than the function $f$ itself. The backpropagation of NeuralODE is based on the adjoint method with a backward ordinary differential equation on the adjoint states $a(t) = \frac{dL}{dh(t)}$. Comparing with traditional DNNs, Neural ODE \cite{Dupont2019-pj} allows continuous time-series modeling and doesn't require discretizing observation and emission intervals.

\subsection{Contributions}
We propose novel learning-based frameworks to solve general optimal control problems. Our frameworks actively incorporates Pontryagin maximum principle (PMP) into the learning process. Our contributions include:
\begin{enumerate}
    \item An optimal planning algorithm (see \cref{discrete_pmp}) with an application to the discrete-time linear quadratic control (LQR) problem with uneven time steps. 
    \item A two-phase-learning framework called \textbf{NeuralPMP} (see \cref{neural_pmp}) that actively uses Pontryagin maximum principle in the training loss function and the learning process. The two phases utilize a variational autoencoder to prevent static learning while enabling the model to learn the correct Hamiltonian and the corresponding dynamics.
\end{enumerate}
To the best of our knowledge, this is the first framework aiming at optimal control problems that utilizes PMP in the learning process instead of merely using PMP to derive closed-form mathematical formulas. 

\subsection{Organization}
The rest of the paper is organized as follows. Section $2$ briefly introduces the mathematical backgrounds including: optimal control formulations, the Pontryagin maximum principle (PMP), and the definitions of reduced Hamiltonian. Section $3$ proposes an algorithm incorporating PMP into the training loss to solve discrete-time control problems. Sections 4 builds on the limitation of such approach in Section $3$ and introduces a two-phase-learning framework called \textbf{NeuralPMP} for continuous-time optimal control. Section $5$ provides experimental results of NeuralPMP on classical control tasks. Finally, section $6$ concludes our work and outlines future research directions.
\section{Optimal control formulation and PMP flow}\label{pmp_theory}
We first state the optimal control problem, and then define the PMP flow that provides the optimal trajectory which, in turn, helps infer the optimal control/action to take.

\subsection{Optimal control formulation}
For simplicity, we work on Euclidean spaces. Our concept can be generalized to a general Banach space. Assume we have state space $X$, a bounded region of $\mathbb{R}^m$, with dual $X^*$, which can be identified with $\mathbb{R}^m$, and a control space $U$, which is another bounded region of $\mathbb{R}^n$. Here $m \in \mathbb{N}$ and $n \in \mathbb{N}$ are the state space and control space dimensions. Given a terminal time $T$, we need to find an optimal control path (also called policy) $u^*$ so that $u^*(t) \in U$ for all time $t \in [0, T]$ and such $u^*$ minimizes the following cost functional:
\begin{equation}
J(q, u) = \int_0^T l(q(t), u(t)) dt + g(q(T))
\end{equation}
Here the function $l = l(q,u)$ is the running cost to move along the control path, and $g = g(q)$ is the terminal cost function. The state $q(t) \in X$ and the control $u(t) \in U$ are subject to following dynamical system constraint:
\begin{align}
\dot{q}(t) &= f(q(t), u(t))\\
q(0) &= q_0
\end{align}
for a fixed starting state $q_0 \in X$. In other words, the optimal control path $u$ satisfied:
\begin{equation}\label{opt_formulation}
u^* = \argmin_{u(t) \in U} J(q, u)
\end{equation}
subject to the given dynamical system constraint. 

\subsection{Dynamic programming principle}
Suppose $\mathcal{V}$ is the value function associated to the control problem \cref{opt_formulation}, and $Q^{s, q}_t$ is the state at time $t$ of the dynamical system starting at $q_s = q$ with the control $\set{u(w)}_{w=s}^t$. The continuous-time dynamic programming for the optimal control problem can be briefly described as:
\begin{equation}
\mathcal{V}(s, q) = \inf_{u: u(w) \in U} \left( \int_s^t l(q(w), u(w)) dw + \mathcal{V}(t, Q^{s, q}_t) \right)
\end{equation}

Notice that the integral term in the continuous-time setting makes it difficult to derive a similar recursive algorithm as the ones in discrete-time value-based reinforcement learning literature. This motivates us to consider the alternative path based on the Pontryagin maximum principle for solving continuous-time optimal control.

\subsection{Pontryagin maximum principle}
As before, let $m$ and $n$ be the state space and control space dimensions. Define the \textit{Hamiltonian} function $H: \mathbb{R}^m \times \mathbb{R}^m \times \mathbb{R}^n \to \mathbb{R}$ as:
\begin{equation}\label{Hamiltonian_def}
H(q, p, u) = p^Tf(q, u) + l(q, u), 
\end{equation}

The Pontryagin maximum principle \cite{Kirk1971-kl} provides the necessary condition for an optimal solution of the general optimal control problem.

\begin{theorem}\label{pmp}(Pontryagin maximum principle) If $(q(.), u(.))$ is an optimal solution of the optimal control problem \cref{opt_formulation}, then there exists the adjoint variable $p(.)$ of $q(.)$ so that:
\begin{align}
&p(T) - \nabla g(q(T)) = 0 \label{pmp1}\\
&\dot{q}(t) = \partial_p H(q(t), p(t), u(t)) \label{pmp2}\\
&\dot{p}(t) = - \partial_q H(q(t), p(t), u(t)) \label{pmp3}\\
&\partial_u H(q(t), p(t), u(t)) = 0 \label{pmp4}
\end{align}
\end{theorem}

The Pontryagin maximum principle (PMP) condition \cref{pmp4} gives us a way to find the optimal control  $u = u(q, p)$ at any specific time, in terms of the generalized coordinates. Now we define the reduced Hamiltonian $h = h(q, p)$ as a function of only the generalized coordinates $(q, p)$ consisting of the state $q$ and the adjoint variable $p$:
\begin{equation}\label{reduced_hamiltonian_def}
h(q,p):= H(q, p, u(q, p))
\end{equation}
where $u$ comes from \cref{pmp4}. Then with the application of the chain rule, PMP conditions \cref{pmp2} and \cref{pmp3} become:
\begin{align}
\dot{q}(t) &= \partial_p h(q(t), p(t))\label{reduced_1}\\
\dot{p}(t) &= -\partial_q h(q(t), p(t))\label{reduced_2}
\end{align}
For instance, $h_q = H_q.1 + H_p.0 + H_u.u_q = H_q + 0 + 0.u_q = H_q = -\dot{p}$. Our main framework in \cref{neural_pmp} focuses on learning the this reduced Hamiltonian function $h = h(q,p)$ and the induced Hamiltonian flow arisen from $h$.
\section{Discrete-time Hamiltonian dynamics learning}\label{discrete_pmp}
Before jumping into our main framework in \cref{neural_pmp}, we propose a discrete-time framework to learn the discretized reduced Hamiltonian dynamics. Moreover, we analyze its advantage over the popular policy-based methods on a specific, irregular linear quadratic control problem.

\subsection{Discrete framework with PMP-based loss}
\textbf{Problem.} Suppose an agent $\mathcal{A}$, initially starts at $q_0$, needs to perform a sequence of actions $u_0, u_1, \cdots, u_{T-1}$ at times $t_0=0, t_1, \cdots, t_{T-1}$ in order to maximize the cost functional $J(q, u)$ introduced in \cref{opt_formulation}.

\textbf{Solution.} From the theory in \cref{pmp_theory}, we build a parameterized neural network $F_{\theta}$ with parameter $\theta$ that maps from initial state $q_0$ to the sequence of optimal actions $(u_0, u_1, u_2, \cdots, u_{T-1})$ by minimizing a Pontryagin maximum principle-based loss function. The training includes 3 steps:
\begin{enumerate}
    \item The initial state $q_0$ is fed into $F_{\theta}$ to obtain a sequence of actions $F_{\theta}(q_0) = (u_0, u_1, u_2, \cdots, u_{T-1})$.
    
    \item The actions $(u_0, u_1, u_2, \cdots, u_{T-1})$ are then fed into a black-box differentiable dynamical model to obtain the  sequence of states $q_1, q_2, \cdots, q_{T-1}, q_T$: $q_k = q_{k-1} + \Delta_k f(q_{k-1}, u_{k-1})$ with the time step $\Delta_k = t_k - t_{k-1}$.
    
    \item The adjoint variables $p_k$ are then calculated backwards in time $p_{k-1}= p_k + \Delta_k H_q(q_k, u_k, p_k)$, where $H$ is the Hamiltonian defined in \cref{Hamiltonian_def}.
\end{enumerate}

\begin{figure}[htbp]
  \centering
  \subfloat[Training loss over 100 iterations \label{fig:training_loss}]{\includegraphics[scale=0.4]{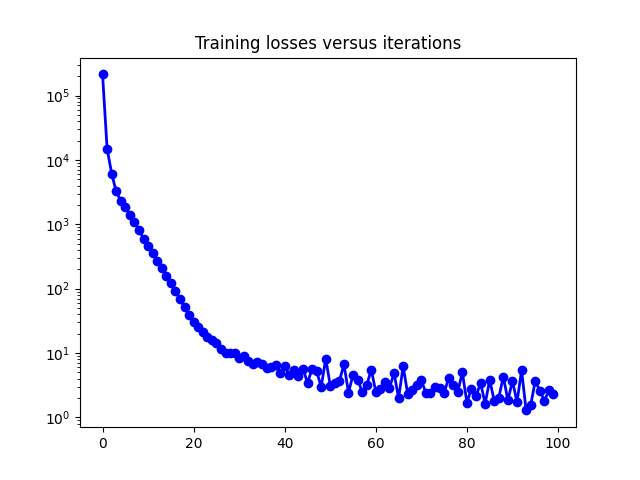}}
  \subfloat[Trained car's behavior \label{fig:car_behavior}]{\includegraphics[scale=0.4]{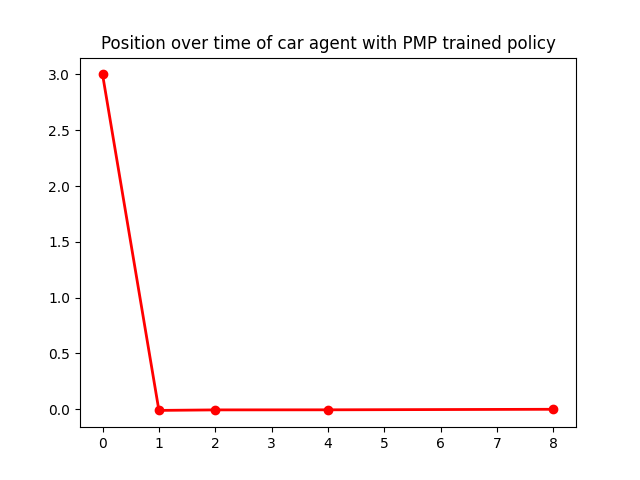}}
\caption{Discrete-time PMP training on LQR problem with uneven timesteps}
\end{figure}

The dynamical system equations \cref{pmp2} and \cref{pmp3} are effectively discretized in the last two steps of the forward pass. We train a network that preserves the last condition \cref{pmp4} of Pontryagin maximum principle (PMP): $\partial_u H(x(t), p(t), u(t)) = 0$. As a result, the loss function for the (policy) network $F_{\theta}$ is defined as:
\begin{equation}
L(\theta) = \sum_{k = 0}^T \norm{H_u(q_k, u_k, p_k)}
\end{equation}

After the training phase, given initial state $q_0$, we produce the optimal behaviors simply by using the sequence of actions: $F_{\theta}(q_0) = (u_0, u_1, \cdots, u_{T-1})$.

\subsection{Special linear quadratic control problem}
We apply the previous model to the following linear quadratic control (LQR) problem:
\begin{equation}
u^* = \argmin J(q, u) \text{, where } J(q,u) = \int_{t_0=0}^{t_T = T} \frac{1}{2}(q^TRq + u^TQu)
\end{equation}
subject to the linear dynamics $\dot{a} = Aa + Bu$. For the special case where state space $X=\mathbb{R}, A = 0, B = 1, R = 1, Q = 0$, the LQR problem is equivalent to the task when a learning agent needs to control its velocity to get as close to the origin as possible even though its initial position can be far from $0$. We train the policy $F_{\theta}$ based on the PMP loss on $100$ iterations with 1600 random samples per iteration. The discrete-time horizon is chosen to be $[0, 1, 2, 4, 8]$. Notice that the times are not equally spaced, and the values of timestamps are taken into the account in this model. The plot (see \cref{fig:training_loss}) shows that the training process converges after the $70$th iteration.

With the actions produced from the trained policy $F_{\theta}$, the trained car agent quickly decreases the velocity in order to get to $0$ as quickly as possible at time $0$. At the subsequent times, the car learns to stay where it was. This is the optimal behavior given that the car needs to minimize its cumulative distance to $0$ over the time period from $0$ to $8$ despite its starting position (see \cref{fig:car_behavior}).

Policy-based algorithms in reinforcement learning learn a policy function $\pi$ mapping from $q_t$ to the optimal action $u_t$. In this special case, we show that no such reasonable function $\pi$ exists. Assume by contradiction, we have such a function $\pi$ so that $u_t = \pi(q_t)$ for each $t$. For $t \geq 1$, $\pi(q_t) = u_t$ must be $0$ because once the car gets to the origin, it should stay there forever. This implies that the policy function is approximately zero function. However, the car must increase or decrease its speed to get to $0$ as quickly as possible right after it starts. Thus, for a starting value $q_0 \neq 0$, $u_0$, which is approximated by $\pi(q_0) = 0$, needs to be non-zero, leading to a contradiction.
\section{Continuous Hamiltonian dynamics learning with forward-backward variational autoencoder}\label{neural_pmp}

The training framework in \cref{discrete_pmp} allows us to gain important insights into the process of actively incorporating PMP principles into the deep learning framework. Nonetheless, the mapping $F_{\theta}$ introduced earlier can be very costly for the time horizon with many steps $t_1, \cdots, t_T$ for large $T$, as it maps to the entire sequence of actions. Moreover, even though the framework can handle uneven time steps, it still requires the discretization of the underlying dynamics, and may not be applicable to a general continuous-time controlled dynamical system.
In this section, by also incorporating PMP conditions into the training process, we build a learning framework for continuous-time optimal control problem with a focus on learning the reduced Hamiltonian function \cref{reduced_hamiltonian_def} and the corresponding Hamiltonian flow.

\subsection{Reduced Hamiltonian learning with PMP-based training loss}
We now describe each components in our learning framework:\\

\begin{figure}[ht]
\begin{center}
\includegraphics[scale=0.55]{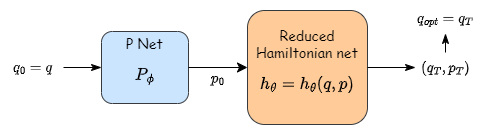}
\end{center}
\caption{Network architecture for phase 1 training}
\label{phase1_architecture}
\end{figure}

\textbf{Neural network architecture:} In this framework, we train two neural networks (see \cref{phase1_architecture}): 
\begin{enumerate}
    \item The parameterized reduced Hamiltonian $h_{\theta} = h_{\theta}(q, p)$ where inputs are the generalized coordinates $(q, p)$.
    \item A $P$-net $P_{\phi}$ that takes input as the starting state $q_0$ and outputs the adjoint variable $p_0$. This $p_0$ is necessary since the Hamiltonian $h_{\theta}$ also needs an adjoint variable as input in addition to the state variable.
\end{enumerate} 

\textbf{Black-boxes environment: } We assume that for the optimal control learning problem, a simulated environment gives us the access to black-boxes evaluations of the dynamics $f$, and/or $f_u$, as well as the cost functions $l$ and $g$. We define $\hat{u}:= \hat{u}(q, p) = -f_u(q)^Tp$ and define $\hat{f}(q_0, \hat{u}_0)$ to be either $f(q_0, \hat{u}_0)$ or $\partial_p h(q(0), p(0))$.\\

\textbf{Training:} The following training procedure is called $\textit{Forward}(T)$ for a fixed terminal time $T$. We first create a replay memory $\mathcal{M}$ that consists of tuples with 5 elements $(q, p, \hat{u}, \hat{f}(q, \hat{u}), l(q, \hat{u}))$. We then train the reduced Hamiltonian net $h_{\theta}$ and the $P$-net $P_{\phi}$ by repeating the following steps:
    \begin{enumerate}
        \item Use NeuralODE \cite{chen2018neuralode} to sample the generalized coordinates $(q_t, p_t)_{t \geq 0}$ from the Hamiltonian flow defined by equations \cref{reduced_1} and \cref{reduced_2} with $h = h_{\theta}$. Then for each generalized coordinates $(q_t, p_t)$ on the Hamiltonian flow, the 5-tuple $(q_t, p_t, \hat{u}(q_t, p_t), \hat{f}(q_t, \hat{u}_t), l(q_t, \hat{u}_t))$ is added to $\mathcal{M}$.
        \item Take a batch of tuples from the memory $\mathcal{M}$, and optimize $h_{\theta}$ with respect to the loss function:
        \begin{align}
            \mathcal{L}^{Forward}(\theta, \phi) &= \alpha_1\norm{P_{\phi}(q_0) - p_0}^2 + \alpha_2\norm{P_{\phi}(q_T)-p_T}^2 \nonumber \\
            &\quad + \beta_1\norm{h_{\theta}(q, p) - (p^T\hat{f}(q, \hat{u}) + l(q, \hat{u}))}^2
        \end{align}
        where $p_T$ is the terminal adjoint variable after running $(q, P_{\phi}(q))$ through Hamiltonian dynamics specified by $h_{\theta}$. The variables $\alpha_1, \alpha_2$ and $\beta_1$ are hyper-parameters and are optimizable by several techniques: Bayesian optimization \cite{bayesian_opt_hyperparam}, Tree-Structured Parzen estimator \cite{tree_parzen_hyperparam}, or bandit algorithms \cite{bandit_hyperparam}. For efficiency and simplicity, we currently employ grid search with early termination. The first two terms are based on the PMP condition \cref{pmp1}, while the last term is used to learn the reduced Hamiltonian correctly.
    \end{enumerate}

\subsection{Second phase of variational training}
In the previous framework, a determistic Hamiltonian flow dictated by the parameterized reduced Hamiltonian $h_{\theta}$ is used. Hence, the induced optimal path is deterministic. We experimentally find that adding noise only reduces the quality of the training and makes it difficult for the agent to learn the correct Hamiltonian. However, to improve the overall policy exploration process, stochasticity needs to be injected into the framework. As a result, in addition to the previous training steps, we propose an addition training phase that makes use of a variational autoencoder that employs forward-backward Hamiltonian dynamics.

\begin{figure}[ht]
\begin{center}
\includegraphics[scale=0.6]{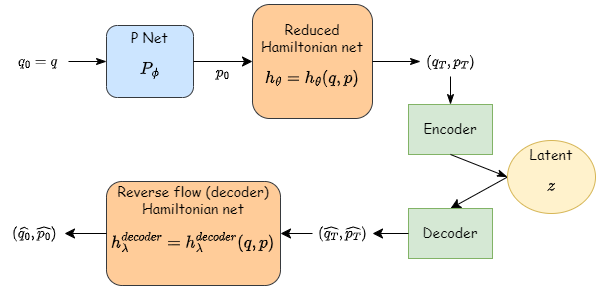}
\end{center}
\caption{Network architecture for phase 2}
\end{figure}

More concretely, after the first training procedure $\textit{Forward}(T)$ is finished, a reverse Hamiltonian flow given by another parameterized reduced Hamiltonian $h^{decoder}_{\lambda} = h^{decoder}_{\lambda}(q, p)$ will be learned. Such learning process include 3 steps:
\begin{enumerate}
    \item We first variationally encode $y:=(q_T, p_T)$ from the $\textit{Forward}(T)$ phase into a latent variable $z \in \mathbb{R}^d$, for a latent dimension $d \in \mathbb{N}$ via an encoder $E_{\lambda_1}$ so that $z = E_{\lambda_1}(y) = E_{\lambda_1}(q_T, p_T)$.
    \item The latent variable $z$ is then mapped to a reconstructed terminal generalized coordinate $\hat{y}:=(\widehat{q_T}, \widehat{p_T})$ via the decoder $D_{\lambda_2}$: $D_{\lambda_2}(z) = (\widehat{q_T}, \widehat{p_T})$.
    \item Using reverse Hamiltonian flow defined by \cref{reduced_1} and \cref{reduced_2} with $h = h^{decoder}_{\lambda}$, we obtain a reconstructed initial generalized coordinate $\hat{x}:=(\widehat{q_0}, \widehat{p_0})$
\end{enumerate}

The training process is similar to that of a variational auto-encoder \cite{Kingma2019-yv} and we optimize the following loss function:
\begin{equation}
\mathcal{L}^{Backward}(\lambda, \lambda_1, \lambda_2) = \norm{x-\hat{x}}^2 + \norm{y-\hat{y}}^2 + \beta_2 \text{KL}(Q(z|y) || P(z))
\end{equation}
The first two terms are reconstruction errors with $y= (q_T, p_T)$, $\hat{y}= (\widehat{q_T}, \widehat{p_T})$, $x = (q_0, p_0)$, $\hat{x}= (\widehat{q_0}, \widehat{p_0})$. The last term is the usual KL-divergence with $P$ and $Q$ being the prior and posterior distributions of the latent variable $z$. The stochasticity now comes from the sampling of $z$ during the training process. $\beta_2$ is a hyperparameter that is again found by grid search with early termination. Note that the reverse flow leading to $(q_0, p_0)$ is based on the \textbf{backward ordinary differential equation} with the dynamics function $(-h^{decoder}_{\lambda})$. We call this additional training phase $\textit{Backward}(T)$ for the terminal time $T$. The final \textbf{NeuralPMP} training procedure includes two phases: \textit{Forward}(T) and then \textit{Backward}(T) trained sequentially with the same terminal time $T$.

\section{Experiments}
We run the continuous-time NeuralPMP framework and compare with other models on 3 classical control tasks in continuous-time settings.

\subsection{Models} 
For benchmarking, we consider 3 different training models:
\begin{enumerate}
\item \textbf{Random Hamiltonian:} This model outputs a random action, which nonetheless follows Hamiltonian dynamics bias rather than totally random.

\item \textbf{NeuralPMP-phase 1}: This model is trained with only the first phase $Forward(T)$ with deterministic reduced Hamiltonian and adjoint neural networks training.

\item \textbf{NeuralPMP:} This is our main framework described in details in the previous \cref{neural_pmp} with 2 phases $Forward(T)$ and $Backward(T)$.
\end{enumerate}

\subsection{Environments:} 
We perform the benchmarking and evaluations on 3 classical control problems with the following evaluation costs:
\begin{enumerate}
    \item \textbf{Mountain car}: The car tries to control the acceleration and subsequently velocity and position to go uphill. The evaluation cost for this environment is the squared distance between the car position and its desired position uphill.
    \item \textbf{Cart pole}: One tries to control the force on the cart to make it move while maintaining balance of the pole on top. The evaluation cost is the pole angle with respect to the vertical direction
    \item \textbf{Pendulum}: One tries to swing the pendulum while maintaining the balance by controlling the torque exerted. The evaluation cost is a normalized form of the pendulum angle.
\end{enumerate}
The smaller the evaluation cost, the better the agents perform. All 3 problems are continuous versions of the environments from the OpenAI Gym suite \cite{openaigym}.

\subsection{Evaluation}
We evaluate 3 models by running each model on $2000$ initial starting points to collect $2000$ trajectories. For each of such trajectories, we calculate the evaluation cost at each time step between $0$ and terminal time $T=1$ with step size $0.05$ (21 time steps in total). For each model and at each time step, we calculate the mean evaluation cost over $2000$ trajectories. We also use bootstrap method to calculate the adjusted variance of the evaluation costs on those $2000$ trajectories by performing 2 steps: 
\begin{itemize}
\item We first collect a bootstrap sample of size $n$ from $2000$ trajectories, where $n$ is a random natural number between $1$ and $2000$. We then calculate the mean on that bootstrap samples. 
\item We perform such a bootstrap calculation 1000 times and calculate the variance of those 1000 mean values obtained.
\end{itemize}

\begin{figure}[htbp]
  \centering
  \subfloat[Mountain Car \label{fig:mountain_car}]{\includegraphics[scale=0.3]{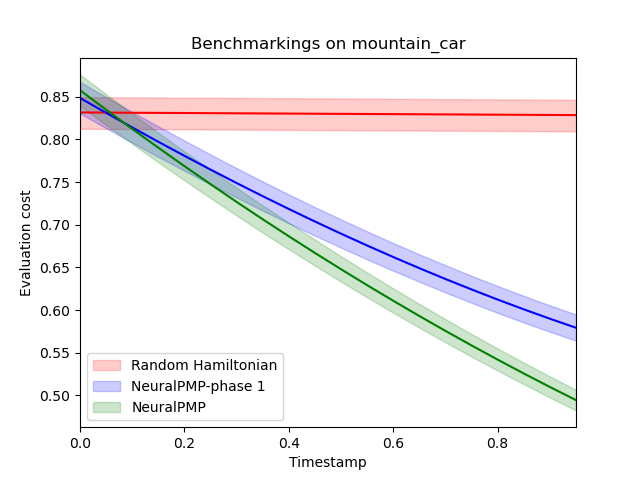}}
  \subfloat[Cart Pole \label{fig:cartpole}]{\includegraphics[scale=0.3]{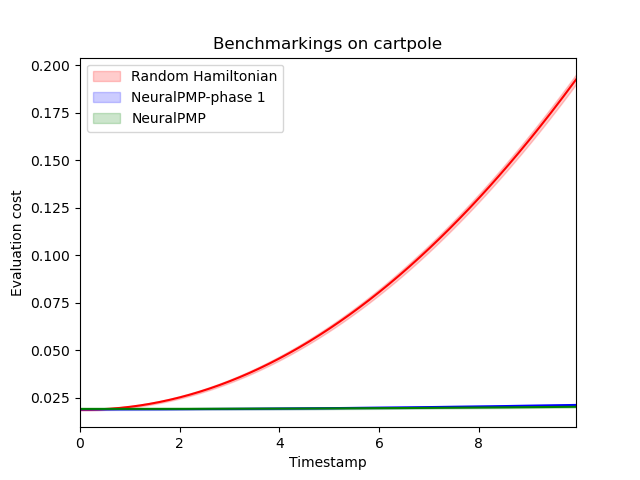}}
  \subfloat[Pendulum \label{fig:pendulum}]{\includegraphics[scale=0.3]{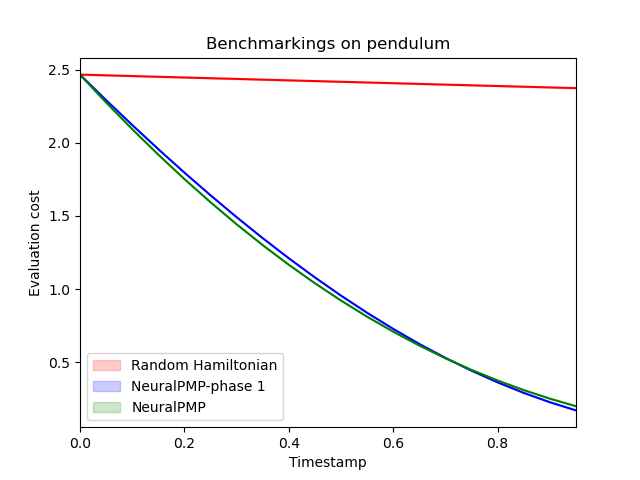}}
\caption{Training benchmarks}
\label{eval_cost_plot}
\end{figure}

The plots (see \cref{eval_cost_plot}) show how evaluation costs (mean and adjusted variance on 2000 trajectories) vary across different time steps for each model on each problem. The plots show that both NeuralPMP models outperform the random Hamiltonian dynamics model. NeuralPMP model performs slightly better than the version with only phase 1 (forward) training. Even in the pendulum case, while the cost is higher at terminal time for NeuralPMP, such cost still drops faster in most time steps than the one from one phase-model. This is due to the stochasticity added via the forward-backward Hamiltonian flows variational auto-encoder structure. This structure overall improves exploration and generalizability while learning the Hamiltonian network accurately. We also report in \cref{metric_table} the average evaluation costs at the terminal time $T=1$ of 3 models on 3 problems. Details about hyper-parameters, neural network architectures, and the continuous-time design of each problem are provided in the appendix.
\begin{table}
    \caption{Average end evaluation costs of different models under 2000 sample trajectories}
    \label{metric_table}
    \centering
        \vspace{3mm}
        \begin{tabular}{c@{\hskip 0.3in}c@{\hskip 0.3in}c@{\hskip 0.3in}c}
        \toprule
        & Random Hamiltonian & NeuralPMP-phase 1 & NeuralPMP \\
        \midrule
        Mountain car & 0.8285 & 0.5793 & \textbf{0.4945}\\
        Cart pole & 0.1924 & 0.0212 & \textbf{0.0202}\\
        Pendulum & 2.375 & \textbf{0.1701} & 0.1986 \\
        \bottomrule
        \end{tabular}
\end{table}
\section{Conclusion}
In this paper, we build a learning framework for optimal control problems by actively using the Pontryagin maximum principle for training and learning and beyond just invocation of a mathematical principle. We first present an optimal planning discrete-time algorithm. We then introduce our main two-phase learning framework called NeuralPMP. This framework allows learning a reduced Hamiltonian dynamics and at the same time derives the optimal control or actions. In the future, we plan to develop model-free versions of our framework to eliminate the need to access black-box dynamics or black-box running/terminal costs. We also would like to extend the framework to larger-scale problems with higher-dimensional states and actions. \\

\textbf{Author acknowledgement:} This research work was supported in part by a grant from the NIH DK129979, in part from the Peter O’Donnell Foundation, the Michael J. Fox Foundation, Jim Holland-Backcountry Foundation, and in part from a grant from the Army Research Office accomplished under Cooperative Agreement Number W911NF-19-2-0333.\\

\textbf{Author Contributions:} This research work is developed by Minh Nguyen and improved by Chandrajit Bajaj. Alphabetical order is currently
used to order author names. The implementation is done by Minh Nguyen and is available at \href{https://github.com/mpnguyen2/neural_pmp}{NeuralPMP project}.

\newpage
\appendix
\section{Network architectures, hyperparameters}
The network architectures for both phases of NeuralPMP are given in \cref{forward_net_architecture} and \cref{backward_net_architecture}. The hyperparameters including $\alpha_1, \alpha_2, \beta_1, \beta_2$ in training loss and number of samples and batch sizes are given in \cref{hyperparamters}.

\begin{table}[ht]
\begin{center}
\begin{tabular}{{c@{\hskip 0.2in}c@{\hskip 0.2in}c@{\hskip 0.2in}c}}
\textbf{Environment}&  \textbf{State dimension} &\textbf{Adjoint network} &\textbf{Hamiltonian network} \\
\hline \\
Cart pole & 4 & [16, 32, 32] & [16, 32, 64, 8]\\
Mountain car & 2 & [8, 16, 32] & [8, 16, 32]\\
Pendulum & 2 & [8, 16, 32] & [8, 16, 32]\\
\end{tabular}
\end{center}
\caption{Network architecture for phase 1 (Forward) training} \label{forward_net_architecture}
\end{table}

\begin{table}[ht]
\begin{center}
\begin{tabular}{{c@{\hskip 0.2in}c@{\hskip 0.2in}c@{\hskip 0.2in}c}}
\textbf{Environment}&  \textbf{Latent dimensions} &\textbf{Hamiltonian decoder} &\textbf{Decoder layers} \\
\hline \\
Cart pole & 4 & [16, 32, 64, 8] & [16, 64]\\
Mountain car & 2  & [8, 16, 32] & [8, 32]\\
Pendulum & 2  & [8, 16, 32] & [8, 32]\\
\end{tabular}
\end{center}
\begin{center}
\begin{tabular}{c@{\hskip 0.3in}c@{\hskip 0.3in}c}
\textbf{Encoder shared layers}&  \textbf{Encoder mean layers} &\textbf{Encoder logvar layers} \\
\hline \\
\text{ }[64] & [16] & [16]\\
\text{ }[32] & [8] & [8]\\
\text{ }[32] & [8] & [8]\\
\end{tabular}
\end{center}
\caption{Network architecture for phase 2 (Backward) training} \label{backward_net_architecture}
\end{table}

\begin{table}[ht]
\begin{center}
\begin{tabular}{c@{\hskip 0.2in}c@{\hskip 0.2in}c@{\hskip 0.2in}c@{\hskip 0.2in}c@{\hskip 0.2in}c@{\hskip 0.2in}c}
\textbf{Environment} & \textbf{$\alpha_1$} & \textbf{$\alpha_2$} &\textbf{$\beta_1$} & \textbf{$\beta_2$} & Number of samples & Batch size\\
\hline \\
Cart pole & 0.1 & 1 & 10 & 1 & 640 & 32\\
Mountain car & 0.1 & 1 & 10 & 1 & 640 & 32\\
Pendulum & 0.1 & 1 & 10 & 1 & 640 & 32\\
\end{tabular}
\end{center}
\caption{Hyperparameters for both phases in NeuralPMP} \label{hyperparamters}
\end{table}

\section{Details of costs and dynamics for control problems}
For the classical control tasks: Cart pole, Mountain car, and Pendulum, we choose the dynamics similar to the one given in the OpenAI Gym suite \cite{openaigym}. The states for these two problems are:
\begin{enumerate}
    \item \textbf{Mountain car}: state $q = (x, \dot{x})$, where $x$ and $\dot{x}$ are position and velocity of the car.
    \item \textbf{Cart pole}: state $q = (x, \dot{x}, \theta, \dot{\theta})$, where $x$ and $\dot{x}$ are position and velocity of the cart, while $\theta$ and $\dot{\theta}$ are angle and angular velocity of the pole with respect to the cart.
    \item \textbf{Pendulum}: state $q = (\theta, \dot{\theta})$, where $\theta$ and $\dot{\theta}$ are angular position and velocity of the pendulum.
\end{enumerate}

The running cost $l(q, u)$ and terminal cost $g(q)$ for 3 problems are:
\begin{enumerate}
    \item \textbf{Mountain car}: $g(q) = (x - x_0)^2 + (\dot{x}-\dot{x}_0)^2$, where $x_0$ and $\dot{x}_0$ are the goal position and velocity. $l(q, u) = 0.1 u^2 + g(q)$ 
    \item \textbf{Cart pole}: $g(q) = \theta^2$, as we don't want the pole to deviate too much from the vertical line. $l(q, u) = 0.5 u^2$.
    \item \textbf{Pendulum}: $g(q) = (\theta+\pi/2)^2$ and $l(q, u) = u^2 + g(q)$.
\end{enumerate}

\newpage
\bibliographystyle{plain}
\bibliography{citation}

\begin{thebibliography}{10}

\bibitem{abbeel2006using}
Pieter Abbeel, Morgan Quigley, and Andrew~Y Ng.
\newblock Using inaccurate models in reinforcement learning.
\newblock In {\em Proceedings of the 23rd international conference on Machine learning}, pages 1--8, 2006.

\bibitem{dpo}
Chandrajit Bajaj and Minh Nguyen.
\newblock {DPO}: Differential reinforcement learning with application to optimal configuration search.
\newblock {\em arXiv preprint arXiv:2404.15617}, 2024.

\bibitem{tree_parzen_hyperparam}
James Bergstra, R\'{e}mi Bardenet, Yoshua Bengio, and Bal\'{a}zs K\'{e}gl.
\newblock Algorithms for hyper-parameter optimization.
\newblock In J.~Shawe-Taylor, R.~Zemel, P.~Bartlett, F.~Pereira, and K.Q. Weinberger, editors, {\em Advances in Neural Information Processing Systems}, volume~24. Curran Associates, Inc., 2011.

\bibitem{AI_pontryagin}
Lucas B{\"o}ttcher, Nino Antulov-Fantulin, and Thomas Asikis.
\newblock {AI} pontryagin or how artificial neural networks learn to control dynamical systems.
\newblock {\em Nat. Commun.}, 13(1):333, 2022.

\bibitem{openaigym}
Greg Brockman, Vicki Cheung, Ludwig Pettersson, Jonas Schneider, John Schulman, Jie Tang, and Wojciech Zaremba.
\newblock Open{AI} gym, 2016.

\bibitem{chen2018neuralode}
Ricky T.~Q. Chen, Yulia Rubanova, Jesse Bettencourt, and David Duvenaud.
\newblock Neural ordinary differential equations.
\newblock {\em Advances in Neural Information Processing Systems}, 2018.

\bibitem{deisenroth2011pilco}
Marc Deisenroth and Carl~E Rasmussen.
\newblock Pilco: A model-based and data-efficient approach to policy search.
\newblock In {\em Proceedings of the 28th International Conference on machine learning (ICML-11)}, pages 465--472. Citeseer, 2011.

\bibitem{Dupont2019-pj}
Emilien Dupont, Arnaud Doucet, and Yee~Whye Teh.
\newblock Augmented neural odes.
\newblock In H.~Wallach, H.~Larochelle, A.~Beygelzimer, F.~d\textquotesingle Alch\'{e}-Buc, E.~Fox, and R.~Garnett, editors, {\em Advances in Neural Information Processing Systems}, volume~32. Curran Associates, Inc., 2019.

\bibitem{PoNNs}
Andrea D’Ambrosio, Enrico Schiassi, Fabio Curti, and Roberto Furfaro.
\newblock Pontryagin neural networks with functional interpolation for optimal intercept problems.
\newblock {\em Mathematics}, 9(9), 2021.

\bibitem{gu2016continuous}
Shixiang Gu, Timothy Lillicrap, Ilya Sutskever, and Sergey Levine.
\newblock Continuous deep {Q}-learning with model-based acceleration.
\newblock In {\em International conference on machine learning}, pages 2829--2838. PMLR, 2016.

\bibitem{heess2015learning}
Nicolas Heess, Gregory Wayne, David Silver, Timothy Lillicrap, Tom Erez, and Yuval Tassa.
\newblock Learning continuous control policies by stochastic value gradients.
\newblock In C.~Cortes, N.~Lawrence, D.~Lee, M.~Sugiyama, and R.~Garnett, editors, {\em Advances in Neural Information Processing Systems}, volume~28. Curran Associates, Inc., 2015.

\bibitem{PDP_end}
Wanxin Jin, Zhaoran Wang, Zhuoran Yang, and Shaoshuai Mou.
\newblock Pontryagin differentiable programming: An end-to-end learning and control framework.
\newblock In H.~Larochelle, M.~Ranzato, R.~Hadsell, M.F. Balcan, and H.~Lin, editors, {\em Advances in Neural Information Processing Systems}, volume~33, pages 7979--7992. Curran Associates, Inc., 2020.

\bibitem{Kingma2019-yv}
Diederik~P Kingma and Max Welling.
\newblock An introduction to variational autoencoders.
\newblock {\em Found. Trends\textregistered{} Mach. Learn.}, 12(4):307--392, 2019.

\bibitem{Kirk1971-kl}
Donald~E Kirk.
\newblock {\em Optimal Control Theory: An Introduction}.
\newblock Prentice-Hall, London, England, 1971.

\bibitem{levine2014learning}
Sergey Levine and Pieter Abbeel.
\newblock Learning neural network policies with guided policy search under unknown dynamics.
\newblock In {\em NIPS}, volume~27, pages 1071--1079. Citeseer, 2014.

\bibitem{bandit_hyperparam}
Lisha Li, Kevin Jamieson, Giulia DeSalvo, Afshin Rostamizadeh, and Ameet Talwalkar.
\newblock Hyperband: a novel bandit-based approach to hyperparameter optimization.
\newblock {\em J. Mach. Learn. Res.}, 18(1):6765–6816, Jan 2017.

\bibitem{Li2017-wy}
Qianxiao Li, Long Chen, Cheng Tai, and Weinan E.
\newblock Maximum principle based algorithms for deep learning.
\newblock {\em Journal of Machine Learning Research}, 18(165):1--29, 2018.

\bibitem{lutter2019deep}
Michael Lutter, Christian Ritter, and Jan Peters.
\newblock Deep lagrangian networks: Using physics as model prior for deep learning.
\newblock In {\em International Conference on Learning Representations}, 2019.

\bibitem{Mnih2016-cy}
Volodymyr Mnih, Adria~Puigdomenech Badia, Mehdi Mirza, Alex Graves, Timothy Lillicrap, Tim Harley, David Silver, and Koray Kavukcuoglu.
\newblock Asynchronous methods for deep reinforcement learning.
\newblock In Maria~Florina Balcan and Kilian~Q. Weinberger, editors, {\em Proceedings of The 33rd International Conference on Machine Learning}, volume~48 of {\em Proceedings of Machine Learning Research}, pages 1928--1937, New York, New York, USA, 20--22 Jun 2016. PMLR.

\bibitem{mnih2015human}
Volodymyr Mnih, Koray Kavukcuoglu, David Silver, Andrei~A Rusu, Joel Veness, Marc~G Bellemare, Alex Graves, Martin Riedmiller, Andreas~K Fidjeland, Georg Ostrovski, et~al.
\newblock Human-level control through deep reinforcement learning.
\newblock {\em nature}, 518(7540):529--533, 2015.

\bibitem{boltyanskiy1962mathematical}
L~S Pontryagin, V~G Boltyanskii, R~V Gamkrelidze, and E~F Mishchenko.
\newblock {\em The mathematical theory of optimal processes, Translated from the Russian by Trirogoff, K.}
\newblock Interscience Publishers John Wiley \& Sons, Inc. New York-London, 1962.

\bibitem{TRPO}
John Schulman, Sergey Levine, Pieter Abbeel, Michael Jordan, and Philipp Moritz.
\newblock Trust region policy optimization.
\newblock In Francis Bach and David Blei, editors, {\em Proceedings of the 32nd International Conference on Machine Learning}, volume~37 of {\em Proceedings of Machine Learning Research}, pages 1889--1897, Lille, France, 07--09 Jul 2015. PMLR.

\bibitem{PPO}
John Schulman, Filip Wolski, Prafulla Dhariwal, Alec Radford, and Oleg Klimov.
\newblock Proximal policy optimization algorithms.
\newblock {\em arXiv preprint arXiv:1707.06347}, 2017.

\bibitem{bayesian_opt_hyperparam}
Jasper Snoek, Hugo Larochelle, and Ryan~P Adams.
\newblock Practical bayesian optimization of machine learning algorithms.
\newblock In F.~Pereira, C.J. Burges, L.~Bottou, and K.Q. Weinberger, editors, {\em Advances in Neural Information Processing Systems}, volume~25. Curran Associates, Inc., 2012.

\bibitem{Sutton1990-qy}
Richard~S Sutton.
\newblock Integrated architectures for learning, planning, and reacting based on approximating dynamic programming.
\newblock In {\em Machine Learning Proceedings 1990}, pages 216--224. Elsevier, 1990.

\bibitem{Todorov2006-jq}
Emanuel Todorov.
\newblock Optimal control theory.
\newblock In {\em Bayesian Brain}, pages 268--298. The MIT Press, 2006.

\bibitem{Wang2019-fk}
Tingwu Wang, Xuchan Bao, Ignasi Clavera, Jerrick Hoang, Yeming Wen, Eric Langlois, Shunshi Zhang, Guodong Zhang, Pieter Abbeel, and Jimmy Ba.
\newblock Benchmarking model-based reinforcement learning.
\newblock {\em arXiv preprint arXiv:1907.02057}, 2019.

\bibitem{zhang2019you}
Dinghuai Zhang, Tianyuan Zhang, Yiping Lu, Zhanxing Zhu, and Bin Dong.
\newblock You only propagate once: Accelerating adversarial training via maximal principle.
\newblock In {\em Advances in Neural Information Processing Systems}, volume~32. Curran Associates, Inc., 2019.

\bibitem{Zhong2019-df}
Yaofeng~Desmond Zhong, Biswadip Dey, and Amit Chakraborty.
\newblock Symplectic {ODE}-{N}et: Learning hamiltonian dynamics with control.
\newblock In {\em International Conference on Learning Representations}, 2020.

\end{thebibliography}

\end{document}